\documentclass{article}
\usepackage{epsfig,array,graphicx}
\usepackage{amsmath,amssymb}
\usepackage{rotating,rotfloat}
\usepackage{pst-sigsys}
\usepackage{pstricks}
\usepackage{multirow}

\title{Two-dimensional generalization of the M\"uller root-finding algorithm and its applications}
\author{Plamen P. Fiziev\thanks{Department of Theoretical Physics, Sofia University ``St. Kliment Ohridski", 5
James Bourchier Blvd., 1164 Sofia, Bulgaria, AND, BLTF, JINR, Dubna, 141980 Moscow Region, Russia,  {\tt fiziev@phys.uni-sofia.bg}, {\tt fiziev@theor.jinr.ru}} and Denitsa R. Staicova\thanks{Department of Theoretical Physics, Sofia University ``St. Kliment Ohridski", 5 James Bourchier Blvd., 1164 Sofia, Bulgaria, {\tt dstaicova@phys.uni-sofia.bg}}}
\begin{document}
\maketitle

\begin{abstract}
We propose a new algorithm for solving a system of two nonlinear transcendental equations with two complex variables based on the M\"uller algorithm. The two-dimensional M\"uller algorithm is tested on systems of different type and is found to work comparably to Newton's method and Broyden's method in many cases. The new algorithm is particularly useful in systems featuring the Heun functions whose complexity may make the already known algorithms not efficient enough or not working at all. In those specific cases, the new algorithm gives distinctly better results than the  other two methods.

As an example for its application in physics, the new algorithm was used to find the quasi-normal modes (QNM) of Schwarzschild black hole described by the Regge-Wheeler equation. The numerical results obtained by our method are compared with the already published QNM frequencies and are found to coincide to a great extent with them. Also discussed are the QNM of the Kerr black hole, described by the Teukolsky Master equation.

\end{abstract}

\pagestyle{myheadings}
\thispagestyle{plain}

\section{Overview}
Solving  a system of two complex-valued nonlinear transcendental equations numerically is a task with varying difficulty, depending on the non-linearity of the system, the types of functions involved and the dimension of the space determined by the system. There are many well-known iterative root-finding algorithms, but most of them are specialized and optimized to work with a narrow set of functions -- for example polynomials or functions with real-valued roots. The two most heavily relied upon one-dimensional algorithms -- the secant method and Newton's method (or the Newton-Raphson method,
\cite{Newton, Raphson, Muller,numerical}) can work with a wide set of complex valued functions under proper conditions (see \cite{numerical}), but they have their weak sides. Newton's method requires the evaluation of the function and its first derivative at each iteration. This increases the computational cost of the algorithm and it makes the algorithm unusable when the procedure evaluating the derivative of the function has numerical problems (for example see the discussion for the Heun functions below) or when derivative becomes zero or changes sign. The secant method avoids this limitation, but in the general case, it has lower order of convergence ($\sim 1.618$) compared to that of Newton's method ($=2$) and the convergence of both of them is strongly dependent on the initial guess.

These problems of the algorithms are inherited by their multi-dimensional versions such as the generalized Newton-Raphson method (\cite{Muller}) and the generalized secant method (Broyden's method, \cite{Broyden}). Although those problems can have varying severity, there are systems in which those algorithms cannot be used effectively. There are also some novel approaches (see \cite{multi}, \cite{cgrasp}), but when they rely on the same one-dimensional algorithms, they are likely to share their weaknesses as well. Clearly there is a need for new algorithms that will enlarge the class of functions we are able to work with efficiently.

A great challenge in front of root-finding algorithms in modern physics can be found in systems including the Heun functions. The Heun functions are unique particular local solutions of a second-order linear ordinary differential equation from the Heun type \cite{heun,heun1_,heun2_,heun3_} which in the general case have 4 regular singular points. Two or more of those regular singularities can coalesce into an irregular singularity leading to confluent differential equations and their solutions: confluent Heun function, biconfluent Heun function, double confluent Heun function and triconfluent Heun function. The Heun functions generalize the hypergeometric function (and also include the Lame function, Mathieu function and the spheroidal wave functions \cite{heun2_,heun3_}) and their wide applications in physics (\cite{heun3_}) was summarized recently in \cite{heun1}. From that paper and the cited therein, it is clear that the Heun functions will be encountered more and more in modern physics from quantum mechanics to astrophysics, hence we need adequate numerical algorithms able to deal with them.

The work with the Heun functions, however, is more than complicated. While there are analytical works on the Heun functions, they were largely neglected until recently and therefore the theory is far from complete. The only software package currently able to work with them is \textsc{maple}. Although, the routines that evaluate them in that package work well in the general case, there are some peculiarities -- there are values of the parameters where the routines break down leading to infinities or to numerical errors. The situation with the derivatives of the Heun  functions is even worst -- in some cases, they simply do not work or their precision is lower than that of the Heun function itself. Also, in some cases there are no convenient power-series representations and then the Heun functions are evaluated in \textsc{maple} using numerical integration. Therefore the procedure goes slowly in the complex domain (compared to the hypergeometric function) which means that the convergence of the root-finding algorithm is essential.

Despite all the challenges in the numerical work with the Heun functions, they offer many opportunities to modern physics. For example, they occur in the problem of quasi-normal modes (QNM) of rotating and non-rotating black holes. In this case, one has to solve a two-dimensional connected spectral problem with two complex equations in each of which one encounters the confluent Heun functions. The analytical theory of the confluent Heun function is more developed than that of the other types of Heun functions, but still many unknowns remain. Again, the evaluation of the derivative of the confluent Heun function is problematic in \textsc{maple} and the behavior of the function is hard to predict. In that situation Newton's method cannot be used as a root-finding algorithm. Broyden's algorithm works well in most cases, but it is slowly convergent even close to a root. It is clear,  then, that  we need a novel algorithm, that will offer quicker convergence than Broyden's algorithm, but without relying on derivatives.

To solve this problem, in the case of a system of two equations in two variables, our team developed a two-dimensional generalization of the M\"uller algorithm.  The one-dimensional M\"uller algorithm (\cite{Muller2}) is a quadratic generalization of the secant method, that works well in the case of a complex function of one variable. It has very good convergence for a large class of functions ($\sim 1.84$) and it is very efficient when the starting point (the initial guess) is close to a root. It is also well convergent when working with special transcendental functions. The two-dimensional M\"uller algorithm seems to inherit some of the advantages of its one-dimensional counterpart like good convergence and usability on large class of functions as our tests show. The new algorithm was used to solve the QNM problem in the case of a Schwarzschild black hole and it proved to work without significant deviations from the results published by Andersson (\cite{Q_N_M}) and Fiziev (\cite{Fiziev1}). Also, preliminary results for the QNM of the Kerr black hole are discussed and for them we also obtain a very good coincidence with published results \cite{special31}.

The article is organized as follows: Section 2 reviews the one-dimensional M\"uller algorithm and its two-dimensional generalization, in Section 3 we discuss some physical application of the method and the numerical results obtained with it and in Section 4 we summarize our results. In the Appendix, the new algorithm is tested on various additional and more simple examples to verify its functionality.

\section{The M\"uller algorithm}
\subsection{One-dimensional M\"uller's algorithm}
The one-dimensional M\"uller algorithm (\cite{Muller,Muller2}) is iterative method which at each step evaluates the function at three points, builds the parabola crossing those points and finds the two points where that parabola crosses the x-axis.  The next iteration is the the point farthest from the initial point.

Explicitly, for every three points $x_{{j-2}}$, $x_{{j-1}}$, $x_{{j}}$ and the corresponding values of the function $f(x)\rightarrow f_{{j-2}}$, $f_{{j-1}}$, $f_{{j}}$, the next iteration $x_{{j+1}}$ is:
\begin{align*}
q={\frac {x_{{j}}-x_{{j-1}}}{x_{{j-1}}-x_{{j-2}}}}, \quad
B=\left( 2\,q+1 \right)f_{{j}}  - \left( 1+q \right) ^{2}f_{{j-1}}+{q}
^{2}f_{{j-2}}\\
C= \left( 1+q \right) f_{{j}}, \quad
{\it D_{1,2}}=B\pm\sqrt {{B}^{2}- 4\,AC}\\
x_{{j+1}}=x_{{j}}-2\,{\frac { \left( x_{{j}}-x_{{j-1}} \right) C}{{
\it D_{max{}}}}},
\end{align*}
\noindent where, $D_{max}$ is the root with bigger absolute value.

The advantages of this algorithm are that it is simple for implementation, it works with complex values, it does not use derivatives and generally, it has higher convergence than the secant method. It also works very well with special functions such as the confluent Heun function. For example the spectra in \cite{Fiziev1} and \cite{spectra} was obtained by the authors using that method.

We will indicate the one-dimensional M\"uller algorithm by the map:
\begin{equation*}
 \mu: \quad \mu(x^{in},F(x))\xrightarrow[{}^{P}]{} x^{fin},
\end{equation*}
\noindent where $x^{in}$, $x^{fin}$ are the starting and end points of the algorithm and the integer $P$ is the number of iterations in which the algorithm completes.

In this notations, the exit-condition for the one-dimensional M\"uller algorithm will be $\mid\!x_{j}-x_{j-1}\!\mid< 10^{-d},$ where $d$ is the number of digits of precision we require.  Combined with checks for the value of the function (if it is small enough), we found this to be the best exit-condition for a root-finding algorithm since it works independently of the actual numerical zero in use, which may vary for the confluent Heun function.

\subsection{Two-dimensional M\"uller's algorithm}
The two-dimensional M\"uller method comes as a natural extension of the one-dimensional M\"uller method.

For two complex-valued functions $F_1(x,y)$ and $F_2(x,y)$ we want to find such pairs of complex numbers $(x_{{}_I},y_{{}_I})$ which are solutions of the system:
\begin{equation}
\begin{cases}
F_1(x_{{}_I},y_{{}_I})=0\\
F_2(x_{{}_I},y_{{}_I})=0
\end{cases}
\label{sys1}
\end{equation}
\noindent where $I=1,\ldots$ numbers the solution in use. From now on, we will omit the index $I$, considering that we work with one arbitrary particular solution. Finding all the solutions of a system is beyond the scope of this article.

Consider the functions $F_1(x,y),\, F_2(x,y)$ as two-dimensional {\em complex} surfaces $z=F_1(x,y)$ and $z=F_2(x,y)$ in a three-dimensional space of the {\em complex} variables $\{x,y,z\}$\footnote{Equivalently, we can consider four {\em real} surfaces in five-dimensional {\em real} hyperspace, which are defined by four {\em real} functions of four {\em real} variables}. Normally, to solve the system, one expresses the relation $y(x)$ from one of the equations, then by substituting it in the other equation, one solves it for $x$ and from $y(x)$ one finds $y$. In the general case, however, this is not possible. The idea of our code is to approximately follow that procedure by finding an approximate linear relation $y(x)$ between the two variables and then using it to find the root of function of one variable trough the one-dimensional M\"uller algorithm.

\begin{figure}
\vspace{-0cm}
\hspace{1.5cm}
\includegraphics[width=300px,height=200px]{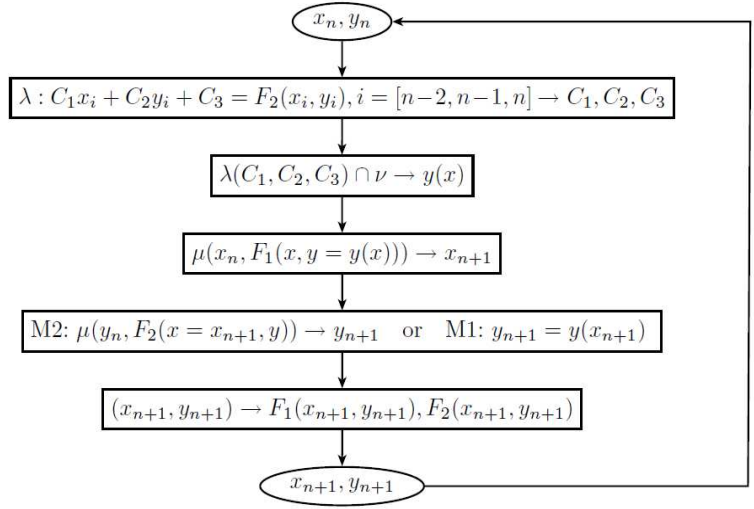}
\caption{A block scheme of the two-dimensional M\"uller algorithm. $\lambda(C_1, C_2, C_3)$ is the plane with equation $z=C_1x + C_2y + C_3$ that crosses trough the 3 pairs of points $(x_{i},y_{i})$ and the function $F_2$ evaluated in them. The plane $\nu$ is defined by the equation $z=0$. The one-dimensional M\"uller algorithm, $\mu(t^{in},F(t)) \rightarrow t^{fin}$, is applied on the function of {\em one variable} $F(t)$ with starting point $t^{in}$ and final point $t^{fin}$}
\label{block}
\end{figure}

To find the linear relation $y(x)$, at each iteration we form the plane passing trough three points of one of the functions and then the equation of the line of intersection between that plane and the plane $z=0$ is used as the approximate relation $y(x)$. This basically means that the so found $y(x)$ is an approximate solution of one of the equations which ideally should be near the real solution in the $z=0$ plane. Substituting this relation in the other function, we run the one-dimensional M\"uller algorithm on it to fix the value of one of the variables, say $x$. Using the value of $x$ in the first function, we again run the one-dimensional M\"uller algorithm on it to fix the value of the other variable -- $y$. Alternatively one can substitute the value of $x$ directly in $y(x)$ to obtain $y$. This ends one iteration of the algorithm.  The process repeats until the difference between two consecutive iterations becomes smaller than certain pre-determined number. This process is systematized on Fig. (\ref{block}).

Explicitly, the code starts by evaluating the two functions $F_{1, 2}(x_i, y_i)$ in three starting pairs of points ($i=1,2,3$) that ideally should be near one of the roots of the system. In our case, those three \emph{initial} pairs are obtained from one starting pair to which we add and subtract certain small complex number. This artificial choice is done only in the first iteration ($n=3$), afterwards we use the output of the last three iterations to form $(x_{n-2},y_{n-2}), (x_{n-1},y_{n-1})$, $(x_n,y_n)$ and the respective $F_{1,2}(x,y)$. Thus on every iteration after $n=3$ the actual complex functions $F_{1, 2}(x_n, y_n)$ are evaluated only once outside of the one-dimensional M\"uler subroutines.

Next we construct the plane passing trough those three points for one of the functions, say $F_2$ by solving the linear system:
\begin{align*}
&C_1x_{n-2}+C_2 y_{n-2}+C_3=F_2(x_{n-2},y_{n-2})\\
&C_1x_{n-1}+C_2 y_{n-1}+C_3=F_2(x_{n-1},y_{n-1})\\
&C_1x_{n}+C_2 y_{n}+C_3=F_2(x_{n},y_{n}).
\end{align*}
\noindent From it one obtains the coefficients $C_1, C_2, C_3$ of the plane $z\!=\!C_1x+C_2y+C_3$.

This plane is intersected with the plane $z=0$ (i.e. $C_1x+C_2y+C_3=0$) and the equation of the line between those two planes is the approximate relation $y(x)$ of the two variables.

We substitute that relation in the first function $F_1(x,y) \rightarrow F_1(x,y(x))$ and we start the one-dimensional M\"uller on that ``linearized`` function of only one variable, $x$. After some pre-determined maximal number of iterations, the exiting point is chosen for $x_{n+1}$ ( $\mu(x_n,F_1(x,y=y(x))) \to x_{n+1}$ \footnote{Since the maximal number of iterations in the one-dimensional M\"uller algorithm is fixed, for simplicity we will omit the index $P$ in this sub-section. The index of the iterations of the two-dimensional M\"uller algorithm is $n$.}) .

Then, there are two possibilities.

Algorithm M1: one could use directly the relation $y(x=x_{n+1})$ to find $y=y_{n+1}$. Or,

Algorithm M2: One can substitute $x=x_{n+1}$ in the other function $F_2(x,y) \to F_2(x=x_{n+1}, y)$ in order to find $y_{n+1}$ using again the one-dimensional M\"uller algorithm ($\mu(y_n,F_2(x_{n+1},y)) \to y_{n+1}$).

Our numerical experiments showed that both approaches lead to convergent procedure.

After $(x_{n+1},y_{n+1})$ are fixed, the two functions $F_{1, 2}(x_{n+1}, y_{n+1})$ are evaluated and if the new points are not roots, the iterations continue.

The exit-strategy in the two-dimensional M\"uller algorithm is as follows:
\begin{enumerate}
 \item To avoid hanging of the algorithm or its deviation to another root we fix maximal number of iterations for the one-dimensional M\"uller subroutine, $P$. Our experience shows that small $P$ (3 --10) is often the more effective strategy, because it stops the one-dimensional M\"uller from straying too far from the actual root of the system.
 \item The precision-condition ($\mid\! x_{j}\!-\!x_{j-1}\!\mid<10^{-d}$) remains in force for the one-dimensional M\"uller. Usually the algorithm exits, because of $j>P$ during the first few iterations of the two-dimensional M\"uller and the closer to the roots it gets, the smaller number of iterations are needed in the one-dimensional M\"uller to reach $d$ and to exit.
 \item The precision $d$ defined by the absolute value of the difference between two consecutive pairs $(x_{n}, y_{n})$ (combined with the values of functions $F_1(x,y), F_2(x,y)$ at them) -- $\mid x_{n} - x_{n-1}\mid < 10 ^{-d}, \mid y_{n} - y_{n-1}\mid < 10 ^{-d}$ is the primary exit-condition of the two-dimensional M\"uller.  When $d$ becomes smaller than certain value, the algorithm exits with a root.
 \item To avoid the two-dimensional M\"uller algorithm from hanging, we set a maximal number of iterations $N$ after which the algorithm reaches exits without fixing a root.
\item A common problem occurs when one of the functions becomes zero before the other function. In those cases, the algorithm accepts the fixed value for a root, say $x^{fin}$, and runs the one-dimensional M\"uller on the other variable until it fixes a root -- $\mu(y_{n},F_2(x^{fin},y)) \to y^{fin}$. The algorithm then exits with a possible root: $(x^{fin},y^{fin})$.
\end{enumerate}
The procedure can be fine-tuned trough change in the starting pair of points, the initial deviation or by switching the places of the functions, or even by replacing the functions with their independent linear combinations.

As we will show in what follows, this method inherits some of the advantages of the one-dimensional M\"uller algorithm, like the quick convergence in proximity of the root and the vast class of functions that it can work with. The major disadvantage comes from the complicated behavior of the two-dimensional complex surfaces defined by the functions $F_{1,2}(x,y)$ which require one to find the best combination of starting points and number of iterations in the one-dimensional M\"uller subroutine so that the algorithm converges to the required root (if it is known or suspected). Generally, it is hard to tell when one point is "close" to a root. In some cases, even if certain starting pair of points is close to a root in terms of some norm, using it as a starting point in the algorithm may still lead to convergence to another root or simply to require more iterations to reach the desired root than if other pair of starting points were used.

It is important to note that unlike Broyden's algorithm and Newton's algorithm which are not dependent on the order of the equations in the system, our two-dimensional M\"uller algorithm depends on the order of the equations. The numerical experiments show that while for some systems, changing the places of the equations has little or no effect on the convergence, in other cases, it slows down or completely breaks down the convergence. While such inherent asymmetry certainly is a weakness of the algorithm, there are ways around it. For example, one may alternate the places of the equations at each iteration or use their independent linear combinations ($F_{1,2}^{*}=\alpha_{1,2} F_1+ \beta_{1,2} F_2$). Those approaches make the algorithm more robust, but since they may cost speed, we prefer to set the order of the equations manually.

A technical disadvantage is that the whole procedure is more CPU-expensive than Newton's method and Broyden's method, since it generally makes more evaluations of the functions -- each one-dimensional M\"uller makes at least $1$ iteration on every step of the two-dimensional M\"uller, thus it makes at least $4$ evaluations of each function. This is because on each iteration of the two-dimensional M\"uller algorithm the functions in use change and thus one cannot use previous evaluations to reduce time. Still, in some cases, as we will show, the so-constructed algorithm is quicker or comparable to Newton's or Broyden's method.

\section{Some applications of the method for systems featuring Heun functions}

We will work only with the confluent Heun functions, which are much better studied than the other types of Heun functions, due to their numerous physical applications. Besides their numerical implementation was used successfully in previous works by the authors. For details on the numerical testing, see the Appendix. 

\subsection{First example}
As a first example, one considers the following system:
\begin{align*}\hspace{+0.1cm}
  F_{1}(x,y)&=\text{HeunC}(\!-\!1.3x,\! 2y,\! 1\!+\!x,\! 4x,\! 1\!-\!1y\!-\!2x^2, \!.75y)\!=\!0\\
  F_{2}(x,y)&=\text{HeunC}(9ix,\! 2.3ix\!+\!y,\! 2ix\!-\!1,\!-\!1.9x(i\!+\!y),\! 2x^2\!+\!2ix\!-\!1.3y\!-\!0.2,\! y)\!=\!0,
\end{align*}
\noindent where HeunC is the confluent Heun function (\cite{heun}) in \textsc{maple} notations.

\begin{table}[ !h]
\centering
\vspace{0cm}
\footnotesize
\addtolength{\tabcolsep}{-1.3pt}{
 \begin{tabular}{|m{1px} |m{53px} | m{122px}  | m{37px}|m{23px} |m{23px} |}
\hline S&  $(x)^{initial}$ & $(x)^{final} $  & $t_{{}_{Broyden}}$[s] & $t_{M2}$[s]& $t_{M1}$[s]\\
  & $(y)^{initial}$ & $(y)^{final} $ & $N_{{}_{Broyden}}$ & $N_{M2}$& $N_{M1}$\\ \cline{1-6}
\hline 
\multicolumn{1}{|m{3px}|}{\multirow{3}{*}{1}}&\multicolumn{1}{c|}{$2.1+0.45i$} &$2.1991016319\!+\!0.2140611770i$&189&161&80\\
&$1.25+0.3i$&$1.2022265008\!+\!0.3588153273i$&14&11(5)*&12(15)*\\\cline{2-6}
\multicolumn{1}{|m{3px}|}{\multirow{3}{*}{}}&\multicolumn{1}{c|}{$2.23+.01i$} &$2.2328663235\!+\!0.0141132493i$&222&91&69\\
&$0.93+0.1i$&$0.9593217208\!+\!0.0508289979i$&23&10(15)*&17(15)*\\
\hline 
\multicolumn{1}{|m{3px}|}{\multirow{3}{*}{2}}&\multicolumn{1}{c|}{$0.49+0.18i$} &$0.4965436315\!+\!0.1849695292i$&208& 102 & 92\\
&$2.001+0.1i$&$1.9999915063\!-\!0.7347653.10^{-5}i$&23&9(5)*&11(4)*\\ \cline{2-6}
\multicolumn{1}{|m{3px}|}{\multirow{3}{*}{}}&\multicolumn{1}{c|}{$0.17+0.97i$} &$0.3495869222\!+\!1.0503235984i$&449&229&244\\
&$2.001+0.1i$&$ 2.0000392386-0.2937407.10^{-4}i$&34&12(5)*&15(5)*\\ \cline{2-6}
\multicolumn{1}{|m{3px}|}{\multirow{3}{*}{}}&\multicolumn{1}{c|}{$0.07+5.147i$} &$0.0608496029\!+\!5.1191008697i$&868&568&489\\
&$2.001+0.051i$&$ 2.0010479243\!-\!0.2491318.10^{-4}i$&36&11(5)*&17(5)*\\
\hline
\end{tabular}}
\caption{$S$ numbers the system in use, $t$ and $N$ label the time and the iterations needed for the algorithms to exit. * denotes the roots dependent on the order of the equations in M1 and M2.\vspace{-0.5cm}}
\label{table0_0}
\end{table}

The results for two pairs of starting points are presented on row S=1 in table \ref{table0_0}. Here, the Newton method do not converge and cannot be used at all. This particularity is due to the behavior of the derivatives of the confluent Heun functions and the known problems of Newton's method near points where the derivative becomes zero or it is not continuous (\cite{numerical}).  From the table, one can see that in this case the M1 modification of the two-dimensional M\"uller algorithm gives the best result, needing less than half the time of the Broyden algorithm to find a root. The M2 modification is slower, but still better than the Broyden algorithm ( $t_{M1}\!<\!t_{M2}\!<\!t_{B}$). 

We continue with two examples from the black hole physics:

\subsection{Physical Application -- QNMs of non-rotating and rotating black holes}

The quasi-normal modes (QNMs) of a black hole (BH) are the complex frequencies ($\omega$) that govern the late-time evolution of the perturbations of the BH metric (\cite{QNM,QNM0,QNM1,QNM2,QNM21}).

\subsubsection*{Second example: Rotating black holes}

To find the QNMs of a rotating black hole, one uses the exact solutions of the Teukolsky radial and angular equations, describing the linearized electromagnetic perturbations of the Kerr metric, in terms of confluent Heun functions, as stated for the first time in full detail in \cite{Fiziev1}. From \cite{Fiziev3}, for the values of the parameters: s=-1, M=1/2, $|r|=110$, m=0, a=0.01, $\theta=\pi/3$, one obtains:

\begin{align*}
\footnotesize
&F_{1}(x,y)=\text{HeunC}(\!-\!1.9996ix, 2.0002ix\!+\!1.0000\!, 0.0002ix\!-\!1.0000\!, 
\!-\!1.9996x(i\!+\!x), \\
&\qquad 1.9995x^2\!+\!1.9998ix\!+\!0.5000\!-\!y\!, \!-\!110.02e^{(4.7124i\!-\!iarg(x))}\!+\!1.0000)
\\ &\qquad \qquad \qquad \times{(110.00e^{(4.7124i\!-\!iarg(x))}})^{(2.00\!+\!0.0002ix)}\!=\!0\\
&F_{2}(x,y)= \frac{\text{HeunC}'(\!0.04x, -1.00, 1.00, -0.04x, 0.50\!-\!1.00y\!+\!0.02x\!-\!0.0001x^2, 0.25)}{\text{HeunC}(\!0.04x, -1.00, 1.00, -0.04x, 0.50\!-\!1.00y\!+\!0.02x\!-\!0.0001x^2, 0.25)}\!+\!\\ &\frac{\text{HeunC}'(\!-\!0.04x, 1.00, -1.00, 0.04x, 0.50\!-\!1.00y\!-\!0.02x\!-\!0.0001x^2, 0.75)}{\text{HeunC}(\!-\!0.04x, 1.00, -1.00, 0.04x, 0.50\!-\!y-0.02x\!-\!0.0001x^2, 0.75)}\!=\!0
 \end{align*}	

For brevity, here the radial equation $F_1(x,y)$ was rounded to only 4 digits of significance. In our numerical experiments, we used the complete system with software floating-point numbers set to $64$, where the derivatives of the confluent Heun functions $\text{HeunC}'$ were replaced with the associate $\delta_N$ confluent Heun function according to equation (3.7) of \cite{Fiziev4}. This was done to avoid the numerical evaluation $\text{HeunC}'$ so that the peculiarities of the numerical implementation of the confluent Heun function (i.e. the use of \textsc{maple} {\em fdiff} procedure) are minimized. The difference in the times needed to fix a root when $\text{HeunC}'$ is used and when it is not used is small for the modes (i.e.$x$) with small imaginary part ($\Delta t\sim15s$), but it increases with the mode number, until it becomes significant for modes with big imaginary part (for the $10^{\text{th}}$ mode -- S=2.3 in table \ref{table0_0} -- the difference is already $\Delta t\sim 100s$). This slowdown is due to the time-consuming numerical integration in the complex domain, needed for the evaluation of $\text{HeunC}'$. 

For that system, three pairs of starting points were used:($0.49+0.18i,2.001+0.1i$), ($0.17+0.97i,2.001+0.1i$), ($0.69e-1+5.146i, 2.001+0.51e-1i$). The results are in row S=2 in table \ref{table0_0}. In this case, once again the two modifications of the two-dimensional M\"uller algorithm M1 and M2 are much quicker than the Broyden algorithm ($t_{M1}\sim t_{M2}<t_{B}$). Newton's method do not converge to a root. The supremacy of the M\"uller algorithms is clear and it is not isolated -- it is observed for other modes or values of the parameters (for example, for $m=1$). As for the precision of the method, the first two modes were found to coincide with at least 9 digits of significance with the already published results of electromagnetic QNMs of a Kerr black hole (see \cite{special31}). We could not find a published value for the third mode. 

\begin{figure}[!h]
\centering
\includegraphics[width=161px,height=150px]{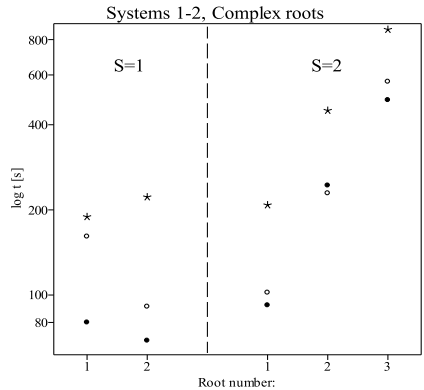}
\caption{A graphical comparison between the times needed by the Broyden method (denoted with asterisks) and the M\"uller methods M1 (the solid circles) and M2 (the circles) }
\label{roots_h}
\end{figure}

Summarizing the results of the previous two examples (rows 1 and 2 in table \ref{table0_0}) in Fig. \ref{roots_h}, one can see that the two-dimensional M\"uller algorithm gives distinctly better times than the Broyden method. Newton's method in those cases cannot be used.

\subsubsection*{Third example: Non-rotating black hole}
Finally, we consider in more detail the problem of the gravitational QNMs of a non-rotating black hole thus testing the new method on a very well studied physical problem. To find the QNMs, one uses the exact solutions of the Regge-Wheeler equations, describing the linearized perturbations of Schwarzschild metric, in terms of confluent Heun functions as stated for the first time in full detail in \cite{Fiziev1}. From \cite{Fiziev1}, when the mass of the BH is set to $2M\!=\!1$, we obtain the following system of the type \eqref{sys1}:
  \begin{align}
F_1&\!=\!(\cos(\theta) \!-\! 1) (\cos(\theta)\! +\! 1) LegendreP(l, 2, \cos(\theta))\!=\!0\label{sys} \\
F_2&\!={\it HeunC} \Big(\!-\!2\,i\omega,2\,i\omega,4,\!-\!2\,{\omega}^{2},4\!-\!l
\!-\!{l}^{2}\!+ \!2\,{\omega}^{2},1\!-\!\mid\! r\!\mid\! e^{-i((\pi\!+\!\epsilon)/2\!+\!\it{arg}(\omega))}\! \Big)\!=\!0 \notag,\hspace{95px}
\end{align}

\noindent where $\omega$ is a complex frequency, $l$ is the angular momentum of the perturbation, $\theta \in [0,\pi]$ is the angle which we set to $\theta=\pi-10^{-7}$ and $|r|=20$. 

The parameter $\epsilon$ is a small variation ($\mid\!\epsilon\!\mid<\!1$) of the phase condition $\arg(r)\!+\!\arg(\omega)\!=\!-\!\pi/2$ (defining the direction of steepest descent, see \cite{Fiziev1}) and we use it to study the stability of the so found QNM frequencies. Since the applicability of the second equation in the system Eqs.(\ref{sys}), $F_2$, may depend on the choice of the values of this parameter, the behavior of the solutions of Eqs.(\ref{sys}) under the variation of $\epsilon$ is still an open problem. Below we present some preliminary results in this direction.

Using Eqs. \eqref{sys}, we run the two-dimensional M\"uller algorithm to find the unknown $l$ and $\omega$. We set the number of digits that \textsc{maple} uses when making calculations with software floating-point numbers to $32$. The precision of the algorithms is still set to $15$ digits.

From the theory, it is known that $l$ is an integer and $l=2, 3 ...$, thus, it offers a way to control the work of the algorithms. We looked for roots only for $l=2$ and we obtained $1.99(9)+1\times10^{-17}i$, with the first different from $9$ digit being the $17^{\text{th}}$.

Because the phase condition  $r\!=\mid\! r\!\mid\! e^{i(-\pi\!+\!\epsilon)/2-i arg(\omega)}$ includes the complex argument in non-analytical way, which cannot be differentiated, this problem cannot be solved directly using Newton's method. Broyden's method, however, also do not work well, because one of the roots $y$ in the pair $(x,y)$ is a real integer, while the other is complex. What happens is that the algorithm fixes the integer root very quickly and the finite differences which are used to evaluate the Jacobian become infinity. Because of this, the algorithm is able to fix only the first $10-11$ digits, while the other algorithms fix $14-15$ digits. Therefore, although Broyden's algorithm gives better times (see Table \ref{table3}) than the two-dimensional M\"uller algorithms due to the less evaluations of the confluent Heun function on each iteration, its precision is much lower and it cannot be increased even by raising the software floating-point number.

\begin{table}[ !htb]
\vspace{-0cm}
\footnotesize
\begin{tabular}{|m{30px}|m{15px}|m{15px}|m{15px}|m{15px}|m{15px}|m{15px}|m{15px}|m{15px}|m{15px}|m{15px}|m{15px}|}
 \hline
Mode:&0&1&2&3&4&5&6&7&8&9&10\\
\hline
$t_B$ [s]&100&99&156&196&386&240&253&282&302&368&398\\
$t_{M2}$ [s]&317&413&595&741&1175&799&874&892&1364&971&1355\\
$t_{M1}$ [s]&202&218&335&357&497&457&396&613&623&594&667\\
\hline
\end{tabular}
\caption{The times needed for Broyden's method ($t_B$) and the two-dimensional M\"uller methods ($t_{M1}$ and $t_{M2}$) to fix a root. Note that while the precision of the former is $10$ digits, the precision of the other two is $14-15$ digits. To obtain those times, we solve the system: $[F_1+F_2,F_1-F_2]$ with starting points: $\omega[n]+0.01+0.01i,2.1+0.01i$, where $n=0..10$. }
\label{table3}
\end{table}

The numerical results are summed in Table (\ref{table2}). In it, the QNM frequencies obtained from Sys. \eqref{sys} are compared to those found by Andersson (\cite{Q_N_M}) with the phase amplitude method. Recently, those results were confirmed by Fiziev (see \cite{Fiziev1}) with the one-dimensional M\"uller method applied on the exact solutions of the radial equation in terms of the confluent Heun function for $l=2$ . To check the accurateness of the new method, we evaluate $\Delta=\mid \omega_{Muller2d}-\omega_{Andersson} \mid$.

\begin{table}[!h]
\vspace{-0cm}
\footnotesize
\begin{tabular}{|m{10px} | m{119px} | m{115px} | m{50px}|}
 \hline n  & Our $\omega$ & Andersson's $\omega$ & $\Delta$\\ \hline
&&&\\
0&0.7473433689+0.177924631i* &0.747343368+0.177924630i&$1.68\times 10^{-9}$\\
1&0.6934219938+0.547829750i* &0.693421994+0.547829714i&$3.60\times 10^{-8}$\\
2&0.6021069092+0.956553966i* &0.602106910+0.956553966i&$1.02\times 10^{-9}$ \\
3&0.5030099241+1.410296405i* &0.503009924+1.410296404i&$1.01\times 10^{-9}$\\
4&0.4150291596+1.893689782i* &0.415029160+1.893689782i&$4.41\times 10^{-10}$\\
5&0.3385988064+2.391216108i &0.338598806+2.391216108i&$9.67\times 10^{-10}$\\
6&0.2665046810+2.895821253i &0.266504680+2.895821252i&$1.48 \times 10^{-9}$\\
7&0.1856446684+3.407682345i &0.185644672+3.407682344i&$3.90\times 10^{-9}$\\
8&-0.030649006+3.996823690i& 0+3.998000i** & 0.0306 \\
9&0.1265270180+4.605289542i&0.126527010+4.605289530i&$1.44\times 10^{-8}$\\
10&0.1531069502+5.121653272i&0.153106926+5.121653234i&$4.52\times 10^{-8}$\\
\hline
\end{tabular}
\caption{A list of the frequencies we obtained for the QNMs of Schwarzschild black hole compared with the numbers found by Andersson. $\Delta\!=\! \mid \!\omega_{Muller2d}-\omega_{Andersson} \!\mid$.  The first 5 frequencies ($n=0-4$, marked with *) were obtained also by Fiziev using the confluent Heun functions and coincide with the presented here except for the $n=1$ where the published by Fiziev value is $0.693421994+0.547829750i$. The 8th mode, marked with **, was obtained by Leaver \cite{Leaver}.}
\label{table2}
\end{table}
One can easily see that most of the modes coincide with precision more than $10^{-8}$. The mode with biggest deviation from the expected value is number 8 in the table \ref{table2}, which is thought to correspond to the so called {\em algebraically special mode}.

Algebraically special (AS) modes have a special place in the QNM studies. The Andersson method is not applicable for them and these are excluded from his consideration. Berti, Cardoso and Starinets (\cite{special3, special1}) make a review on the results so far concerning these modes. Theoretically the 8th mode should be purely imaginary with value $4i$, if it indeed corresponds to the AS case. Leaver (\cite{Leaver}) found that mode at frequency $0.000 000 + 3.998 000i$. In our results the 8th mode has small but non-zero real part. This clearly shows that mode does not agree with the theoretical predictions for the AS mode (see also \cite{special2}).

Furthermore, the mode with $n=8$ depends critically on the value of $\epsilon$ in the phase condition. For  $\epsilon \in [-0.1,0.1]$, our method could not find any corresponding frequency. Outside this interval, there is  such mode, but the sign of its real part depends on the sign of $\epsilon$ -- positive epsilon leads to positive real part and vice versa ($\omega_{n=8}\!=\!\text{sgn}(\epsilon)\,0.030649006+3.996823690i$).

We examine more closely the relation $\omega_n(\epsilon)$ for all modes:

For modes $n\leq4$ there is no dependence on $\epsilon$ and they come in pairs symmetrical to the imaginary axis\footnote{The ranges where each mode is found, however, depend on $\epsilon$: for $n=0$: $\epsilon \in [\mp 0.8, \pm 0.75]$, for $ n=1: \epsilon \in [\mp 0.8, \pm 0.45]$, for $ n=2: \epsilon \in [\mp 0.8, \pm 0.25]$, for $ n=3:\epsilon \in [\mp 0.8, \pm 0.1]$, for  $ n=4: \epsilon \in [\mp 0.8, 0]$, where the first sign corresponds to frequencies with positive real part, the second sign -- to those with negative real parts. The imaginary parts for each mode $n$ coincide.}: $\omega_n=\pm|Re(\omega_n)|+Im(\omega_n)i$ . These results confirm the roots for $n=0,1,2,3,4$ published in \cite{Fiziev1}.

Modes with $n>4$ (but $n\neq8$) depend on $\epsilon$ and the symmetry is broken -- there is only one frequency from each mode corresponding to certain $\epsilon$, i.e. $\omega_n(\epsilon)=-\text{sgn}(\epsilon)|Re(\omega_n)|+Im(\omega_n)i$. The dependence of $|Re(\omega_n)|$ and $Im(\omega_n)$ on $\epsilon$ requires more careful investigation and is beyond the scope of this work. Outside the above mentioned ranges, the modes $n>4$ disappear or get translated.	

The so found relation $\omega_n(\epsilon)$ may be a signal for some kind of instability that needs to be studied more carefully. The presented here results are preliminary and the detailed study will be published elsewhere. Here our goal is just to illustrate the ability of our new numerical algorithm to tread this problem. For the case $n=8$, similar behavior was mentioned by other authors, too \cite{special3, special1,special2}.

From the presented results, it is clear that the two-dimensional M\"uller algorithms work very well for the complicated system of equations Eqs.(\ref{sys}) and it even allowed the study of the dependency $\omega_n(\epsilon)$ done here for the first time. The next step is to apply the method to the more complicated physical problem of QNM of rotating black hole (see Example 2.). As recently discussed by Fiziev (see \cite{Fiziev3}), rotating black holes are described by the radial and angular Teukolsky equations which can be solved
analytically in terms of confluent Heun functions. Preliminary results show that our algorithm works well in this situation too.

\section{Conclusion}
We presented the general idea of a method for solving a system of two complex-valued nonlinear transcendental equations with complex roots based on the one-dimensional M\"uller method. As it was shown, the new method inherits some of the advantages of its one-dimensional counterpart as good convergence for a large class of functions. The efficiency of the new method is illustrated on different examples (For examples including elementary functions and simple special functions see Appendix.).  The numerical results showed that in many cases, it is comparable to Newton's method and Broyden's method, while in some specific systems  -- like those featuring Heun functions -- it is quicker than the other two. The complete mathematical investigation of the proposed new method, and especially its theoretical order of convergence under proper conditions on the class of functions $F_1,F_2$ is still an open problem.

Furthermore, the two-dimensional M\"uller algorithm works well with special functions like the confluent Heun function in terms of which one can solve analytically the Regge-Wheeler equation \cite{Fiziev1}, the Zerilli  equation \cite{Fiziev2}, the Teukolsky radial and angular equations \cite{Fiziev3} and many others. Using this algorithm, we were able to solve directly the problem of quasi-normal modes of a Schwarzschild black hole with higher precision than that of the Broyden method. The so found solutions agree to great extent with previous published numerical results thus confirming the usefulness of the method. The new algorithm was used to study the stability of the so found QNM with respect to small variations in the phase condition in the radial variable and some preliminary results were presented. For other applications of the method see the recent references \cite{PP, PP1}.

Future application of the two-dimensional M\"uller method would be in the case of QNMs of a Kerr black hole described by the radial and angular Teukolsky equations, which will be published elsewhere.

\section*{Acknowledgements}
The authors are grateful to prof. Hans Petter Langtangen for the critical reading of the early version of the manuscript and for the useful advices.

This article was supported by the Foundation "Theoretical and
Computational Physics and Astrophysics", by the Bulgarian National Scientific Fund
under contracts DO-1-872, DO-1-895, DO-02-136, and Sofia University Scientific Fund, contract 185/26.04.2010.

\section*{Author Contributions}
P.F. gave the idea and the outline of the two-dimensional generalization of the M\"uller algorithm,
chose the physical problem on which to test the algorithm and he supervised the project.

D.S. is responsible for the realization of the algorithm in
\textsc{maple} code, the testing and the optimization of the code and for the numerical results presented here.

Both authors discussed the results of the tests of the algorithm, commented them and were involved
in trouble-shooting of the code at all stages. The manuscript was prepared by D.S. and edited by P.F..

\section*{Appendix}
\section{Numerical testing}
All the algorithms are realized as procedures on the software package  \textsc{maple}, the tests are done on  \textsc{maple} 15, on Linux x64, CPU Intel Centrino Core 2 Duo, on 2.2GHz.  The number of digits that \textsc{maple} uses when making calculations with software floating-point numbers is set to $64$. For Newton's method and Broyden's method we used the analytical formulas \cite{Muller}, where the Jacobian in both cases is evaluated exactly or with finite differences respectively (i.e. without the Sherman-Morrison formula).

The times presented below are obtained after running each procedure $10$ times using the garbage collection function $gc();$ in \textsc{maple} on each calculation, so that each run represents an independent numerical experiment. The total time for each method is then divided by $10$ and rounded to 3 digits of significance. This way, even though the times depend on the system load at the moment, they are representative for the four methods in each example. The notable exception of this ''averaging`` are all the systems featuring Heun functions, where such procedure would require too much time and thus they are evaluated only once, using the function $gc();$. 

The precision in all the examples is $15$ digits, but only the first $10$ after the decimal point are presented here. In all the cases, the initial deviation where needed is $0.001$. Some of the examples are from \cite{num_an} p.617-618.

The numerical results for the test-systems are summarized in Table \ref{table0_1} in the Appendix. In it, we compare Newton's method, Broyden's method and the two versions of the two-dimensional M\"uller algorithm discussed in section 2.2. -- $M2$ which uses \emph{two} one-dimensional M\"uller subroutines to fix the $(x_{n+1},y_{n+1})$ and $M1$ which uses \emph{one} one-dimensional M\"uller subroutine to find $x_{n+1}$ and then it evaluates directly $y_{n+1}=(-C_3-C_1x_{n+1})/C_2$. The number in the brackets in the $N_{M1}$ and $N_{M2}$ columns is $P$, the maximal iterations in the one-dimensional M\"uller subroutine. Everywhere in the table, for each $(x^{in},y^{in})$, the four algorithms exit with the same $(x^{fin},y^{fin})$ with precision of 15 digits.

\subsection{Elementary functions}
\begin{enumerate}
 \item[1.] $F_{1}(x,y)=y^2+3x-5+x^2=0$\\
$F_{2}(x,y)=x^2+3y-1=0$\\
  \item[2.] $F_{1}(x,y)=x(1-x)+4y=12$\\
$F_{2}(x,y)=(x-2)^2+(2y-3)^2=25$\\
\end{enumerate}
For those systems (rows S=1,2 in table \ref{table0_1}),	 the number of iterations and the time needed to find a root in the two-dimensional M\"uller algorithms (for M1 in particular) are generally close to those of Broyden's algorithm ($t_{N}\!<\!t_{B}\!\approx \!t_{M1}\!<\!t_{M2}$). For real roots, however, the algorithms M1 and M2 are the {\em quickest} of the four.

\subsection{Trigonometric, exponential and logarithmic functions}

\begin{enumerate}  
\item[3.] $F_{1}(x,y)=y-1/4\sin(x)-1/4\cos(y)=0$\\
$F_{2}(x,y)=5x^2-y^2=0$\\
  \item[4.] 	$F_{1}(x,y)=exp(-3x)\cos(y)+x=0$\\
$F_{2}(x,y)=x^2-3yx+y^2=0$\\
\end{enumerate}
Again, for real roots, the two-dimensional M\"uller methods are quicker than both other methods (rows S=3,4 in table \ref{table0_1}). Newton's method is much quicker when the initial conditions and the roots are complex.

It is important to discuss the dependency of two-dimensional M\"uller methods (M1 and M2) from the maximal number of iterations in the one-dimensional M\"uller subroutine, $P$. In many cases, changing $P$ only affects the time needed for the algorithms to complete. There are cases, however, where this parameter becomes critical. For example, when one uses as starting points (4.4-5.0i, 8.5-16i) on the system $S=5$, the following 7 roots are obtained:
\footnotesize
\begin{align*}
&r_0\!=\!(0.3487096094\!+\!0.4633971546i, 0.9129336096\!+\!1.2131895010i)\\
&r_1\!=\!(3.0248444374\!-\!4.3689275542i, 7.9191455477\!-\!11.4380008313i)\\
&r_2^{\pm}\!=\!(0.1632674377\!\pm\!0.6065137375i, 0.0623626119\!\pm\!0.2316676331i)\\
&r_3\!=\!(1.1119158619\!-\!1.8296636950i, 2.9110335191\!-\!4.7901217415i)\\
&r_4\!=\!(4.0158133827\!-\!5.6039287836i, 10.5135359284\!-\!14.6712760260i)\\
&r_5\!=\!(-5.3999170768\!-\!3.12\times10^{-37}i, \!-\!14.1371664435\!-\!4.38\times10^{-43}i)
\end{align*}
\normalsize
From them, Newton's method converges to the root $r_0$ (after 25 iterations), Broyden's method -- to $r_1$ (after 27 iterations). With the two-dimensional M\"uller methods (M1 and M2) depending on the parameter $P$ one obtains:
\begin{itemize}
 \item with M2: $r_2^{+}$ for $P=3, 9, 10, 12-14$ and $P>16$, $r_0$ for $P=5$, $r_5$ for $P=6$, and $r_2^{-}$ for $P=8$ (after averagely 10 iterations),
 \item with M1: $r_3$ for $P=3$, $r_4$ for $P=4-10$ and $P>16$ and $r_1 $ for $P=11-14$ (after averagely 11 iterations).
\end{itemize}

Similar behavior is observed in the next example:
\begin{enumerate} 
\item[5.] 	$F_{1}(x,y)=\ln(x^2+y^2)-\sin(yx)-\ln(2)+\ln(\pi)=0$\\
$F_{2}(x,y)=e^{x-y}	+\cos(yx)=0$.\\
\end{enumerate}
For real starting points, Newton's method is not convergent, since it remains locked to the real axis. Broyden's method also did not exit with a root for real starting points. The two-dimensional M\"uller methods on the other side, when started from (0.5,0.5) gave the root $r_1^-$ (see below). Other purely real initial conditions either gave a root or the algorithm did not converge.

Using the complex initial conditions (2.27+0.001i, 1.27), the following four roots were found:
%
\footnotesize
\begin{align*}
&r_0\!=(0.2129109625\!-\!2.4380400935i, \!-\!1.3216238026\!-\!4.6551486236i)\\
&r_1^{\pm}\!= (0.9203224533\!\pm\!0.7487874838i, 1.4188731053\!\mp\!0.5453380689i)\\
&r_2\!=(\!-\!1.6645201248\!+\!1.380553001i, 1.66452012482\!+\!1.38055300197i)
\end{align*}
\normalsize

\normalsize Newton's method exits with $r_0$ after 24 iterations. Broyden's algorithm is not convergent. From the two-dimensional  M\"uller algorithms one obtains:
\begin{itemize}
 \item  with M2 after averagely 12 iterations: $r_1^{-}$ for $P=3$, $r_2$ for $P=4$ and $r_1^{+}$ for $P\ge5$ .
 \item  with M1 after averagely 14 iterations: $r_1^{+}$ for $P=3$ and $P\ge5$ and $r_1^{-}$ for $P=4$.
\end{itemize}

From the last two systems it is clear that $P$ represents an additional parameter of the two-dimensional M\"uller algorithm. It can be used to improve convergence, but in some cases, it can lead to different roots for the same initial conditions. Such instability depends on the system and it can be avoided by starting the algorithm closer to the root.	

\subsection{Special functions}
Finally, we consider the following two systems:
\begin{enumerate}
\item[6.]  $F_{1}(x,y)=x^2-y+5\sin(x-2)=0$\\
$F_{2}(x,y)=\text{J}(3, y)+5x-3=0$\\
\item[7.]  $F_{1}(x,y)=x^7-e^y+{{}_2F_2}([1], [3], x^2-3x)=0$\\
$F_{2}(x,y)={H^1}(7, y+1-x)=0$\\
\end{enumerate}
\noindent where $J()$ is the Bessel functions of the first kind, $H^{1}()$ is the Hankel function of the first kind and ${}_2F_2$ is the generalized hypergeometric function.

In this case (rows S=6,7 in table \ref{table0_1}) the two-dimensional M\"uller algorithms are comparable to Newton's algorithm, while Broyden's algorithm is often the quickest of the four. This is likely due to the computational burden of the derivative or of each additional function evaluation. Note, however, that our goal is not to have an algorithm that is better than Newton's method, but to have an algorithm that has good convergence and that does not need to evaluate derivatives. In that, the performance of the new algorithms is satisfactory, especially since in some cases like 6.2 and 6.3, M\"uller's algorithms are the quickest.
\subsection{Discussion}

\begin{figure}[htbp]
\centering
$\begin{array}{ccc}
\hspace{-0.3cm}\includegraphics[width=161px,height=150px]{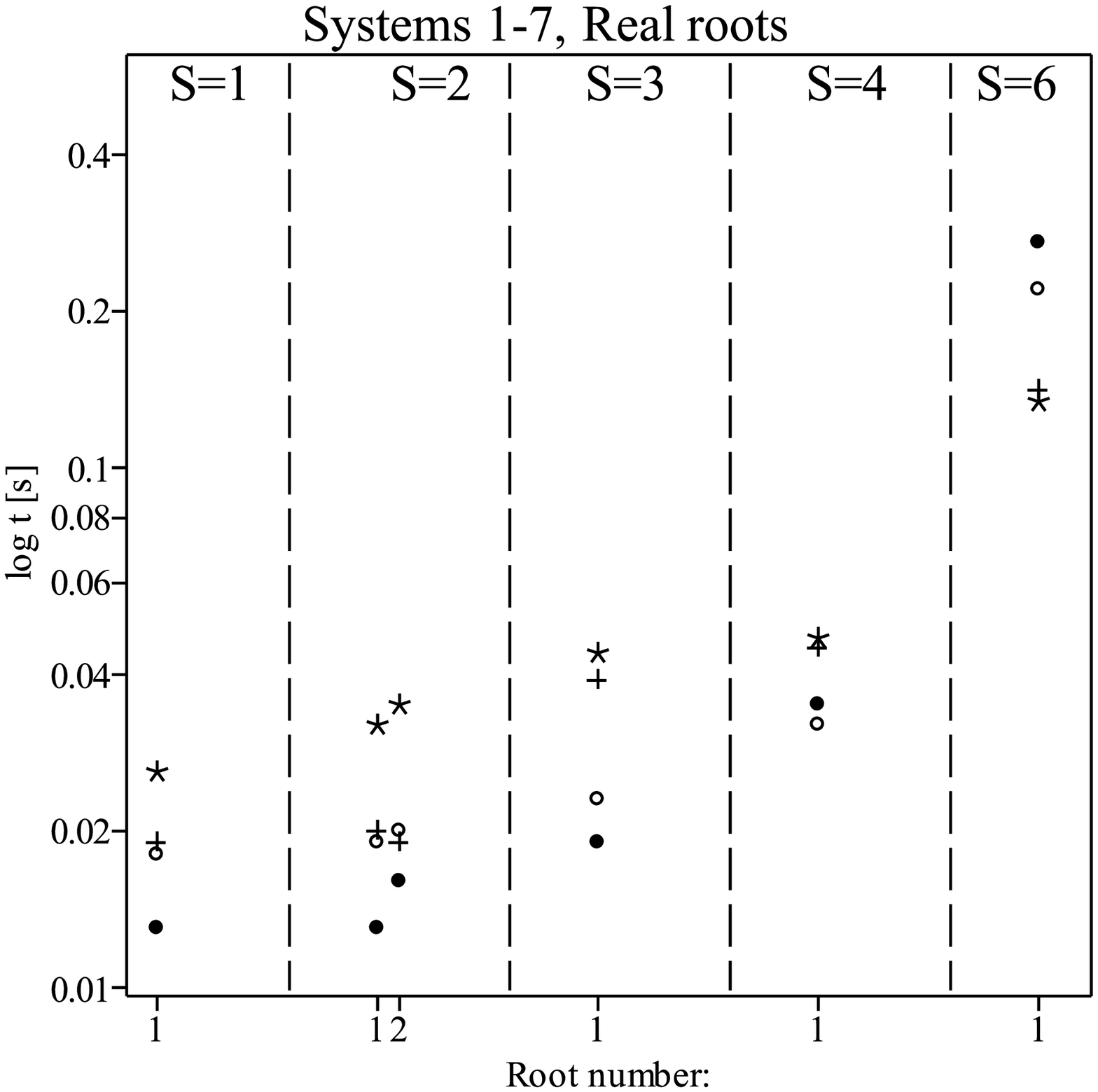}
&\includegraphics[width=161px,height=150px]{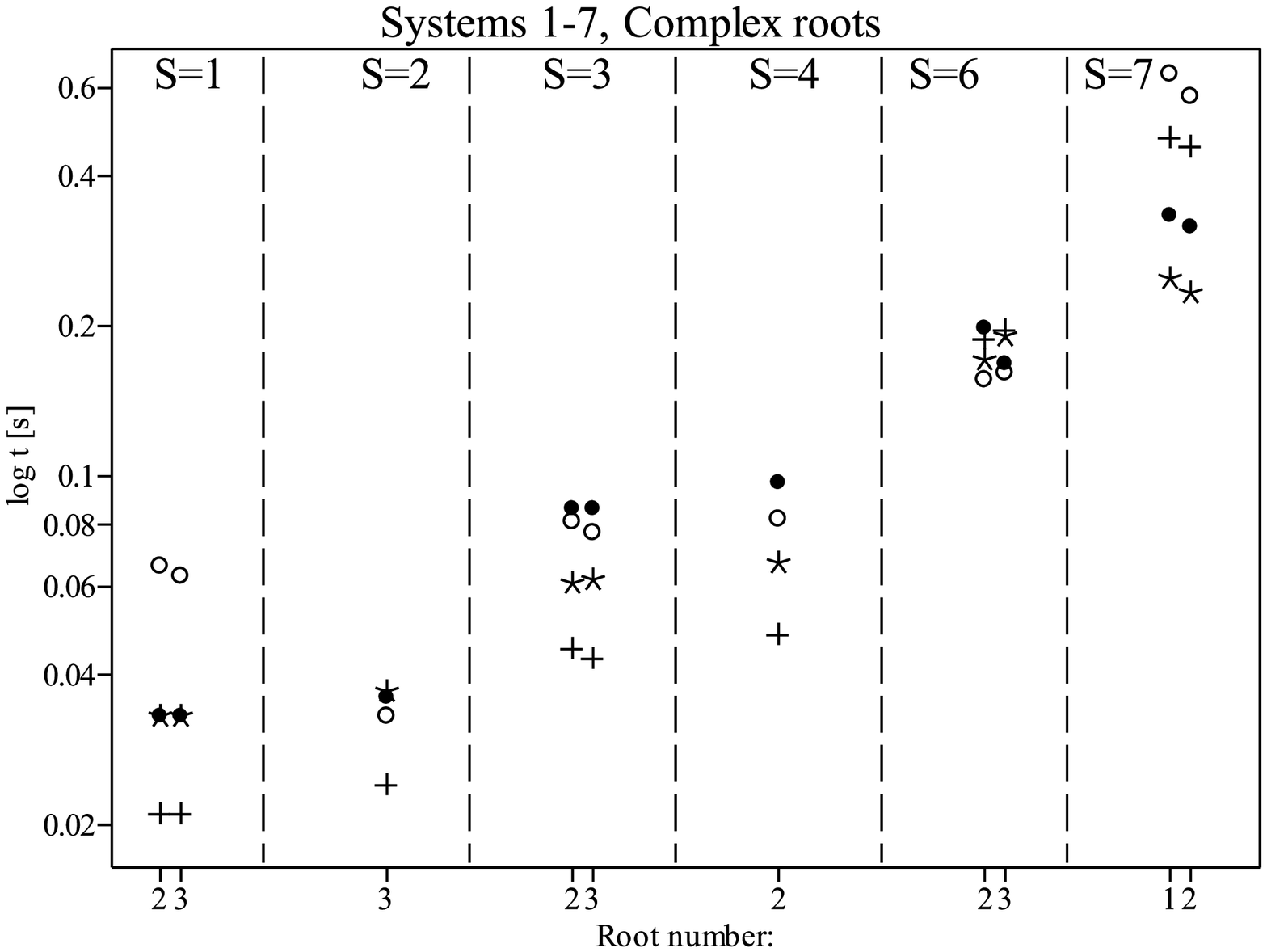}
\end{array}$
\caption{A graphical comparison between the times needed by the Newton method, the Broyden method and the M\"uller methods M1 and M2. With crosses we denote Newton's method, with asterisks -- Broyden's method, with solid circles -- M1, with circles -- M2 }
\label{roots}
\end{figure}

The numerical investigations above (see Fig.\ref{roots} and  also table \ref{table0_1}) 
show that in general, the two modifications of the two-dimensional M\"uler algorithm work comparably well to the more established algorithms -- Newton's and Broyden's, even if sometimes they require more time and iterations to fix a root. In some specific cases, like those with real roots or those featuring confluent Heun functions (see Section 3), however, the two-dimensional M\"uler algorithms are often the quickest of the four. 

It was demonstrated that while in Newton's and Broyden's algorithms the exit points depend only on the starting points (when the initial deviation is fixed), the two-dimensional M\"uller algorithm depends also on the number of iterations in the one-dimensional M\"uller subroutine ($P$). Surprisingly, the time for fixing a root do not depend in a straight-forward way from $P$, since sometimes increasing $P$ leads to decreasing of the total time. This can be expected, because $P$ is the maximal number of iterations in the one-dimensional M\"uller subroutine, but the actual number of iterations depends on the precision. Therefore, this is one more way to fine-tune the algorithm to achieve a known or suspected root.

The examples also showed the importance of the order of the functions in the system (see the systems in the table \ref{table0_1} marked with * and $\dagger$). While in most cases all the four algorithms find at least one root of the system, there are initial conditions which lead to a divergence or to ''undesired`` root. This problem can be avoided by starting the procedure closer to the root, changing the places of the equations in the system or using a linear combination of the functions (say $F_{1}\! \to\! F_1\!+\! F_2, \, F_{2} \!\to\! F_1\!-\! F_2$). 
We conclude that even though the new algorithms admit some further improvements and developments, they work well enough to be tested in real physics problem.

\begin{table}[ !h]
\vspace{0cm}
\footnotesize
\addtolength{\tabcolsep}{-1.3pt}{
 \begin{tabular}{|m{1px} |m{53px} | m{122px} | m{35px} | m{37px}|m{23px} |m{23px} |}
\hline S&  $(x)^{initial}$ & $(x)^{final} $ & $t_{{}_{Newton}}$[s] & $t_{{}_{Broyden}}$[s] & $t_{M2}$[s]& $t_{M1}$[s]\\
  & $(y)^{initial}$ & $(y)^{final} $ & $N_{{}_{Newton}}$ & $N_{{}_{Broyden}}$ & $N_{M2}$& $N_{M1}$\\ \cline{1-7}

\multicolumn{1}{|m{1px}|}{\multirow{3}{*}{1}}&\multicolumn{1}{l|}{$1.689$} &$1.1890465736$&0.019&0.026&0.018&0.013\\
 &$-0.637$ &$ -0.1379439181$&8&10&10*(3)&8(3)\\ \cline{2-7}

 \multicolumn{1}{|m{1px}|}{}&\multicolumn{1}{l|}{$1.321\!+\!3.520i$}&$0.8214691720\!+\!3.5201983985i$& 0.021&0.033&0.066&0.033\\
 &$3.738\!-\!1.927i$&$4.2389950548\!-\!11.9278229759i$&8&10&12(3)&8(3)\\
\cline{2-7}

 \multicolumn{1}{|m{1px}|}{}&\multicolumn{1}{l|}{$1.321\!-\!3.520i$}&$0.8214691720\!-\!3.5201983985i$& 0.021&0.033&0.063&0.033\\
 &$3.738\!+\!1.927i$&$4.2389950548\!+\!1.9278229759i$&8&10&12(3)&8(3)\\
\hline
\multicolumn{1}{|m{1px}|}{\multirow{2}{*}{2}}&\multicolumn{1}{l|}{$-.5$} &$-1.0000000000$&0.020&0.032&0.019&0.013\\
 &$3$   &$ 3.5000000000$&8&10&9(4)&9(3)\\ \cline{2-7}
\multicolumn{1}{|m{1px}|}{}&\multicolumn{1}{l|}{$3.046$}&$2.5469464699$& 0.019&0.035&0.020&0.016\\
&$3.484$&$3.9849974627$&8&10&10(3)&9(3)\\ \cline{2-7}
\multicolumn{1}{|m{1px}|}{}&\multicolumn{1}{l|}{$0.726\!+\!4.335i$}&$0.2265267650\!+\!4.3352949767i$& 0.024&0.037&0.033&0.036\\
&-$2.242\!-\!0.592i$&-$1.7424987313\!-\!0.5927935709i$&8&9&7(6)&8(6)\\
\hline
\multicolumn{1}{|m{3px}|}{\multirow{3}{*}{3}}&\multicolumn{1}{l|}{$0.621$} &  $0.1212419114$&0.039&0.044&0.023&0.019\\
&  $-0.228$ & $0.2711051557$&10$\dagger$&13$\dagger$&9(3)$\dagger$&8(4)*$\dagger$\\ \cline{2-7}
\multicolumn{1}{|m{3px}|}{}&\multicolumn{1}{l|}{-$0.422\!+\!1.476i$}&-$.9222203725\!+\!1.4764038337i$& 0.045&0.061&0.081&0.086\\
&-$2.562\!+\!3.301i$&-$2.062147443\!+\!3.3013393343i$ &9&11&8(4)&11*(3)\\\cline{2-7}
\multicolumn{1}{|m{3px}|}{}&\multicolumn{1}{l|}{ $1.468\!-\!1.635i$}&$0.9685241736\!-\!1.6351708695i$& 0.043&0.062&0.077&0.086\\
&-$2.665\!+\!3.656i$&-$2.1656858901\!+\!3.6563532190i$&9&11&7(5)&11*(3)\\
\hline
\multicolumn{1}{|m{3px}|}{\multirow{3}{*}{4}}&\multicolumn{1}{l|}{$-.35$} &$-0.5600551872$&0.045&0.047&0.032&0.035\\
&$-1.05$&$-1.4662435158$&11&12&7(4)&10(4)\\ \cline{2-7}
\multicolumn{1}{|m{3px}|}{}&\multicolumn{1}{l|}{$0.55-0.6i$}&$0.3487096094\!-\!0.4633971546i$& 0.048&0.067&0.082&0.097\\
&$1.14-1i$&$0.9129336096\!-\!1.213189501i$&10&12&7(6)&13*(3)\\
\hline
\multicolumn{1}{|m{3px}|}{\multirow{3}{*}{6}}&\multicolumn{1}{c|}{$1.2+0.09i$} &  $.6863031247$&0.141&0.134&0.220&0.271\\
&$-5.5+0.01i$ &   -$4.3646459533$&9$\dagger$&11$\dagger$&9*(4)$\dagger$&11(3)$\dagger$\\ \cline{2-7}
\multicolumn{1}{|m{3px}|}{}&\multicolumn{1}{c|}{$7.2-3.6i$}&$5.8404591703\!-\!3.0854927956i$& 0.188&0.171&0.156&0.198\\
&-$ 11.9+5.001i$&-$10.6712592035\!+\!5.7445552813i$ &11&14&10*(3)&14*(4)\\	\cline{2-7}
\multicolumn{1}{|m{3px}|}{}&\multicolumn{1}{c|}{-$5.1-1.006i$}&-$4.9297777922\!-\!1.1922443124i$& 0.196&0.191&0.161&0.168\\
& $ 16.0+5.51i$&$17.4620338366\!+\!5.7870418188i$&11&15&11(3)&13(3)\\
\hline
\multicolumn{1}{|m{3px}|}{\multirow{3}{*}{7}}&\multicolumn{1}{c|}{$1.1-.45i$} &$0.8288091244\!-\!0.4046494664i$&0.476&0.249&0.640&0.333\\
&-$2.4-4.2i$&-$2.3507488745\!-\!4.6830120304i$&10&13&11(3)&12(3)\\\cline{2-7}
\multicolumn{1}{|m{3px}|}{}&\multicolumn{1}{c|}{$.5-.87i$}&     $0.2656154750\!-\!0.8757700972i$& 0.457&0.233&0.577&0.316\\
&-$3.21-5.14i$&-$2.9139425238\!-\!5.1541326612i$&10&13&8*(4)&11*(3)\\%
\hline
\end{tabular}}
\caption{$S$ numbers the system in use, $t$ and $N$ label the time and the iterations needed for the algorithms to exit. * denotes the roots dependent on the order of the equations in M1 and M2. In the $\dagger$ case, the places of the equations were switched to obtain that root.\vspace{-0.5cm}}
\label{table0_1}
\end{table}

\newpage

\end{document}